\documentclass[b5paper, 12pt] {article} 
\usepackage{hilbert10} 
\begin{document}

\title{A positive solution to Hilbert's 10th problem} 
    
\author{Michael Pfender\footnote{michael.pfender@alumni.tu-berlin.de}}
\date{April 2014\\last revised \today}
\maketitle

\abstract{
Polynome codes and code evaluation; arithmetical theory frames;
$\mu$-recursive race for decision; decision correctness;
decision termination in Zermelo-Fraenkel set theory ZFC+ with
axiom of choice and consistency provability; decision correctness 
in theory $\T = \PR$ of Primitive Recursion; 
comparison with the negative result of Matiyasevich; positive solution
for each single diophantine polynomial in 
\emph{p.\,r.\ non-infinite-descent theory} $\piR=\PR+(\pi).$
}

\tableofcontents

\section*{Introduction}

Within theory $\ZFC^+ = \ZFC+\Con_{\ZFC}$ of Zermelo-Fraenkel 
\textbf{set} theory with axiom of choice $\AC,$ strengthened by 
formula $\Con_{\ZFC}$ which is to express $\ZFC$'s \emph{internal,} 
\emph{g\"odelised} consistency, we solve Hilbert's 10th
problem positively: we organise decision of
diophantine polynome codes---decision on overall 
\emph{non-nullity}---as an enumerative $\mu$-recursive race for 
a (first) zero (\emph{counterexample}), against race for a first 
internal $\ZFC$-non-nullity \emph{proof} for a given such polynomial
code, given as the (nested) list of coefficients. Comparison with
Matiyasevich's negative solution of Hilbert's 10th problem gives
inconsistency of theory $\ZFC+\Con_{\ZFC}$ whence self-inconsistency
$\ZFC \derives \neg \Con_{\ZFC}.$ 

In a final section we
plug our positive solution of the problem into the constructive
framework of \emph{p.\,r.\ non-infinite descent theory} $\piR=\PR+(\pi)$
out of \emph{Arithmetical Foundations} in the References.

This is to give a decision algorithm for each single 
diophantine equation (in a uniform way), as asked in the original 
Hilbert's 10th problem.

\section{Hilbert's 10th Problem}

We attempt a positive solution to Hilbert's 10th problem. In its original 
form it reads:

\bigskip
10. DETERMINATION OF THE SOLVABILITY OF A DIOPHANTINE EQUATION
Given a diophantine equation with any number of unknown quantities
and with rational integer numerical coefficients: \emph{To devise a process
according to which it can be determined by a finite number of operations
whether the equation is solvable in rational integers.} 

[translation quoted from \NAME{Matiyasevich} 1993.]

\bigskip
Formally, this text allows for a separate decision algorithm 
(``process'') for each diophantine polynomial. But it is clear that
a decision-\emph{family} must be \emph{uniform} in a suitable sense.

\emph{Correctness} of our alleged $\mu$-recursive decision algorithm 
$\nabla_{\ZFC}: \mr{DIO} \parto \two=\set{0,1}$ builds, within $\ZFC^+,$ 
on diophantine soundness inferred by $\Con_{\ZFC}$ over $\ZFC.$ 
Termination follows from (countable) Choice. This already within 
$\ZFC.$ Together this gives the wanted decision $\nabla=\nabla_{\ZFC}$ 
within $\ZFC^+,$ of all polynome codes in $\mr{DIO}\subset\N.$

Comparison with Matiyasevich's negative Theorem,
\emph{unsolving} Hilbert's 10th Problem, 
theorem in particular of (classically quantified 
Arithmetical Theory) $\ZFC^+,$ gives a contradiction within 
$\ZFC^+,$ hence {self-inconsistency} of $\ZFC,$ 
and from that in particular $\omega$-inconsistency.

In a final section we show correctness and irrefutable termination
of \emph{localised} decision $\nabla[D]$---for each single diophantine
polynomial $D=D(\vec{x})$---within the constructive framework of
p.\,r.\ \emph{finite-descent-theory} $\piR=\piR+\Con_{\piR}$ out
of op.\,cit. 

\section{Polynome coding and code evaluation}

Diophantine polynomials $D = D(\vec{\xi}): \Z^* \to \Z$ (``in $\DIO$'') 
are \LaTeX/\textbf{ASCII} coded into  
  $$\mr{DIO} \defeq \Z^{\an{*}} \iso \union_{m\geq1} \Z[\xi_1,\ldots,\xi_m]  
                 = \union_{m \geq 1} \Z[\xi_1][\xi_2] \ldots [\xi_m]$$ 
as nested coefficient \emph{lists} $\Z^{\an{*}} \subset \N.$

\smallskip
\noindent
$[\,$The symbols $\xi_i$ are the \emph{indeterminates.}$\,]$

\subsubsection*{Example:}
\begin{align*}
D = D(\xi_1,\xi_2)  
& = (2\cdot\xi_1^{\,0}+3\cdot \xi_1^{\,1}- 4\cdot\xi_1^{\,3})\cdot\xi_2^{\,0} \\ 
& \quad
    +(0\cdot\xi_1^{\,0}+3\cdot\xi_1^{\,1}-7\cdot\xi_1^{\,2})\cdot\xi_2^{\,1}  
                                           + (1-4\cdot\xi_1)\cdot\xi_2^{\,2}
\end{align*}
is coded 1-1 as (nested) {coefficient list} 
\begin{align*}
\ccode{D}  
& = \an{\an{2;3;0;4};\an{0;3;-7};\an{0};\an{1;-4}}: \\ 
& \one \to \mr{DIO} \bydefeq \Z^{\an{*}} \subset \N: \\
& \text{\emph{defined element, point} of $\mr{DIO}$}
\end{align*}

\subsubsection*{PR evaluation of $\mbf{DIO}$ codes:}

Evaluation $\mr{ev} = \mr{ev}(d,\vec{x}): \mr{DIO} \times \Z^*$
is PR {defined} 
\begin{align*}
& \mr{ev}(d,\an{\vec{x};x_{m+1}\,}) = \mr{ev}(d,\an{x_1;\ldots;x_m;x_{m+1}\,})\\
& \defeq \mr{ev}(\horner(d,x_{m+1}),\an{\,\vec{x}\,}): \\
& \mr{DIO} \times \Z^* \supset \Z[\vec\xi,\xi_{m+1}] \times \Z^{m+1} 
    \ovs{\iso} (\Z[\vec\xi] [\xi_{m+1}] \times (\Z^m \times \Z) \\
& \xto{\iso} (\Z[\vec\xi][\xi_{m+1}] \times \Z) \times \Z^m 
              \ovs{\horner \times \id} \Z[\vec\xi] \times \Z^m \xto{\ev} \Z,
\end{align*} 
recursively by iterative application of Horner's schema
to the hitherto trailing argument, until all of the arguments 
(constants or variables) are substituted into 
their corresponding indeterminates $\xi_j.$

\smallskip
Result then is the integer $\mr{ev}(d,\vec{x}),$
constant or integer variable.

\smallskip
For the \textbf{example} above, $D = D(\xi_1,\xi_2),$ with argument string  
$\an{x_1;x_2} :\,= \an{23;64} \in \Z^*,$ we get
\begin{align*}
& \mr{ev}(d,\an{x_1;x_2}) 
    = \mr{ev}(\an{\an{2;3;0;4};\an{0;3;-7};\an{0};\an{1;-4}},\an{23;64}) \\
& = \horner(\,
    ((((-4\cdot64+1)\cdot\xi_1+0))
                       \cdot\xi_1 +(-7\cdot64+3)\cdot64)\cdot\xi_1 \\
& \qquad\qquad\qquad\qquad\qquad\qquad\quad
    + ((4\cdot64+0)\cdot64+3)\cdot64+2\,,23\,) \\
& = ((((-4\cdot64+1)\cdot23+0))\cdot23 +(-7\cdot64+3)\cdot64)\cdot23 \\
& \qquad\qquad\qquad\qquad\qquad\qquad
    + ((4\cdot64+0)\cdot64+3)\cdot64+2
\end{align*} 
\textbf{First step:} apply Horner's schema to coefficient list 
$d \in \mr{DIO}$ und (trailing) Argument $x_2:$ {indeterminate}
$\xi_1$ is {coded} by {list nesting} and is seen 
as a \emph{constant}, as an element of intermediate ring $\Z[\xi_1]:$ 
  $$\Z[\xi_1,\xi_2] \bydefeq \Z[\xi_1][\xi_2] \bydefeq (\Z[\xi_1])\,[\xi_2].$$  
\textbf{Last---here second---step:} evaluation of 
$\Z[\xi_1]$ polynomial in remaining indeterminate $\xi_1$ on 
remaining argument $x_1,$ by a last application of Horner's schema.  

\section{Arithmetical frame theories}

We consider here as frame theories---for our decision algorithm --
\textbf{on one hand} classically quantified arithmetical theories 
$\T = \bfQ+\AC$ with (countable) axiom of choice, as in particular 
Zermelo-Fraenkel set theory $\T = \ZFC = \ZF+\AC.$ Frame then is the strengthening  
  $$\T^+ = \T+\Con_\T = \ZFC+\Con_{\ZFC}$$
of $\T$ by its own consistency-\emph{formula} 
\begin{eqnarray*}
\Con_\T &=& \neg\,(\exists\,k \in \N)\,\Pro_\T(k,\code{\false}) \\
&=& (\forall k)\,\neg\,\Pro_\T(k,\code{\false})\ (\text{G\"odel}),
\end{eqnarray*}
see \NAME{Smorynski} 1977 and op.\,cit.

Strengthening by this consistency formula will provide for 
\emph{correctness} of our \emph{decision process} (Hilbert).

\smallskip
\textbf{On the other hand} we take as frame the 
Free-Variables (categorical) theory $\T = \PR = \PRa$ of 
\emph{Primitive Recursion with predicate abstraction into subsets} 
$$(\chi = \chi(a): A \to \two) 
  \boldsymbol\mapsto \set{A:\chi}=\set{a\in A:\chi(a)}$$
out of op.\,cit.\,,\ 
$\T = \bfS$
in Smorynski's notation, as well as \emph{descent theory} 
$\piR = \piR^+ = \piR+\Con_{\piR}:$ that theory is self-consistent,
$\piR \derives \Con_{\piR,}$ main result of op.\,cit.

\section{A $\mu$-recursive race for decision}

We \textbf{define} an enumerative \emph{race}---for
$d \in \mr{DIO}$ thought \emph{passive, fixed,} and $k \in \N$ 
\emph{running}---for satisfaction of 
\begin{align*}
& \ph_0(d,k) = [\,\ev(d,\ct_*k)=0\,]\ \text{against} \\
& \ph_1(d,k) = \Pro_{\T}(k,\code{(\vec{x})\ev(d,\vec{x})\neq 0}):
  \mr{DIO} \times \N \to \two = \set{0,1}, \\
& \ct_* = \ct_*\,k: \N \ovs{\iso} \Zlist
  \ \text{Cantor-type \emph{count},}\  
    \vec{x} \in \Z^*\ \text{free under code.}
\end{align*}

This race towards \emph{termination} is defined as 
a---formally partial---$\mu$-recursive mapping as follows
within the theory $\hatT$ of \emph{partial PR maps,} i.\,e.\ of
(partially defined) \emph{$\mu$-recursive maps,} cf.\ again op.\,cit.:
  $$t = t(d) = \mu\set{k\,|\,\ph_0(d,k)\,\lor\,\ph_1(d,k)}: 
                                          \mr{DIO} \parto \N. \quad (*)$$

\bigskip
\textbf{Decision candidate} then is
\begin{align*}
\nabla d =
& \begin{cases}
  0\ \myif\ \ph_0(d,t(d)) \\
  1\ \myif\ \ph_1(d,t(d)) \\
  \end{cases} \\
& \quad \\
=
& \begin{cases}
  0\ \myif\ \ev(d,\ct_*(t(d)))=0 \\
        \qquad (\emph{zero found}) \\
  1\ \myif\ \Pro_{\T}(t(d),\code{\ev(d,\vec{x})\neq 0}) \\
        \qquad (\text{\emph{internal proof} found
                        for \emph{global non nullity}})
  \end{cases} \\
& : \mr{DIO} \overset{(\id,t)} {\parto} \mr{DIO} \times \N \to \two.
\end{align*}

\textbf{Question:} Is $\nabla$ \emph{well-defined} as a partial map?
In which frame?

\subsubsection*{Well-definedness of the decision within 
$\T^+ = \ZFC^+=\ZFC+\Con_{\ZFC} = \T+\Con_\T:$} 
\begin{align*}
\T^+\ \derives\ 
& \ph_0(d,k)\,\land\,\ph_1(d,k') \\
&     \qquad (\text{\emph{cases-overlap} 
                        {$\mathit{Assumption}$}}) \\
& \implies \ev(d,\ct_*k)=0 \\
& \quad\,\land\,\Pro_\T(k',\code{(\vec{x})\,\ev(d,\vec{x})\neq 0}) \\
& \implies \Pro_\T(j(k,k'),\code{\false}) \\
& \implies \neg\,\Con_\T \implies \false, \\
& j=j(k,k'): \N^2\to \N\ \text{suitable.}
\end{align*}

\textbf{Consequence:}  
  $$\T^+\ \derives\ \neg\,[\ph_0(d,k)\,\land\,\ph_1(d,k')\,]: 
                                     \mr{DIO} \times \N^2 \to \two,$$
$\nabla = \nabla_{\T}(d): \mr{DIO} \parto \N$ is \emph{well-defined} 
as a (\emph{formally partial}) $\mu$-recursive map, within 
$\T^+ = \T+\Con_\T.$ 

\subsubsection*{Well-definedness of decision within \emph{descent} 
theory $\piR:$}

We consider now \emph{descent theory} $\piR$ out of op.\,cit.\ 
strengthening $\PR$ by axiom $(\pi)$ of 
\emph{non-infinite endo driven descending complexity with complexity
values in polynomial semiring $\N[\omega],$} and its logical properties,
in particular \emph{soundness} giving $\piR\derives \Con_{\piR}.$
 
Decision $\nabla=\nabla_{\piR}(d):\mr{DIO}\parto \two$ 
is in fact well-defined as a partial PR map, within theory $\piR,$ 
since---in parallel to the above case $\T = \ZFC:$
\begin{align*}
\piR\ \derives\ 
& \ph_0(d,k)\,\land\,\ph_1(d,k') \\
&     \qquad (\text{\emph{cases-overlap} 
                        {$\mathit{Assumption}$}}) \\
& \implies \ev(d,\ct_*k)=0 \\
& \quad\,\land\,\Pro_{\piR}(k',\code{(\vec{x})\,\ev(d,\vec{x})\neq 0}) \\
& \implies \Pro_{\piR}(j(k,k'),\code{\false}) \\
& \implies \text{``$\neg\,\Con_{\piR}$''} \implies \false, \\
& j=j(k,k'): \N^2\to \N\ \text{suitable.}
\end{align*}
The latter since $\piR\derives \Con_{\piR}.$  

\subsubsection*{Well-definedness of DIO-decision within $\PR$ itself}
Decision $\nabla=\nabla_{\PR}(d):\mr{DIO}\parto \two$ 
is well-defined as a partial PR map, within theory $\hatPRa$ of
partial PR maps since
\begin{align*}
\hatPRa\ \derives\ 
& \ph_0(d,k)\,\land\,\ph_1^{\DIO}(d,k') \\
&     \qquad (\text{\emph{cases-overlap} 
                        {$\mathit{Assumption}$}}) \\
& \Iff \ev(d,\ct_*k)=0 \\
& \quad \land\,\Pro_{\DIO}(k',\code{(\vec{x})\,\ev(d,\vec{x})\neq 0}) \\
& \implies \Pro_{\DIO}(j(k,k'),\code{\false}) \\
& \implies \false, \\
& j=j(k,k'): \N^2\to \N\ \text{suitable.}
\end{align*}
The latter by \emph{diophantine soundness} of $\T=\PR,$ see 
\NAME{Smorynski} 1977, \textsc{Theorem} \textbf{4.1.4}.

\section{Decision Correctness}

\textbf{Decision Correctness, result-0-case:}
\begin{align*}
\T \derives\ 
& [\,\ph_0(d,t(d)) 
    \implies \mr{ev}(d,\ct_*\,\circ\,t(d))=0\,] \\
& \subseteq\,\true_{\mr{DIO}}: 
    \mr{DIO} \overset{(\id,t)} {\parto} \mr{DIO} \times \N \to \two:
\end{align*}
\textbf{If} race-for-decision $\nabla$ \emph{terminates} on 
DIO-code $d,$ with \textbf{result} $0,$ \textbf{then} 
(evaluation of) $d$ has (at least) one zero, namely 
$$\ct_*\,\circ\,t(d) \in \N.$$ 

\subsubsection*{Correctness, result-1-case:}
\begin{align*}
\T \derives\ 
& \ph_1(d,k) \implies \Pro_{\DIO}(k,\code{\ev(d,\vec{x})\neq 0}) \\
& \implies \ev(d,\vec{x})\neq 0: 
  (\mr{DIO} \times \N) \times \Z^* \to \two, \\
& (d \in \mr{DIO},\ k \in \N,\ \vec{x} \in \Z^*\ \text{all free}), \\
& \quad
    \text{or, with quantifier decoration:} \\
\T \derives\
& (\forall\,d \in \mr{DIO})(\forall\,k \in \N)
    (\forall\, \vec{x} \in \Z^*) \\
& [\,\ph_1^{\T} (d,k) \implies \Pro_{\DIO}(k,\code{\ev(d,\vec{x})\neq 0}) \\
& \implies \ev(d,\vec{x})\neq 0\,].
\end{align*}
\textbf{If} race-for-decision $\nabla$ \emph{terminates} on 
DIO-code $d,$ with \textbf{result} $1,$ \textbf{then} 
(evaluation of) $d$ has no zeroes. 

This because of \emph{Diophantine Soundness} of $\T,$ see 
\NAME{Smorynski} 1977, \textsc{Theorem} \textbf{4.1.4} again. 

\subsubsection*{Correctness in result-1-case, 
under termination condition:}
Substitution of $t(d)$ for $k$ in the above gives 
\begin{align*}
\T^+,\piR,\PR\ \derives\ 
& [\,\ph_1^{\DIO} (d,t) \implies \ev(d,\vec{x})\neq 0\,] 
  \subseteq\,\true_{\mr{DIO} \times \Z^*},\\
& d \in \mr{DIO},\ \vec{x} \in \Z^*\ \text{both free}: 
\end{align*}
\emph{Correctness} of $\nabla(d)$ where defined, in \emph{both} 
defined cases: in case of reaching
\textbf{result} 0, as well as in case of reaching  
\textbf{result} 1.

\smallskip
[\,For {partial maps} $f,\,g: A \parto B,$ 
$f\,\subseteq\,g$ designates \emph{inclusion} of the
\emph{graphs} of $f$ and $g.$]


\section{Termination}

We show first

\smallskip
\textbf{Pointwise non-derivability of non-termination:}

For no diophantine \emph{point} $d_0: \one \to \mr{DIO}$ 
$\T$ derives non-termination of $t$ at $d_0.$

\textbf{Proof:}
\inference{}
{ {$\mr{Assumption}$} \\
& $\T \derives\ (\vec{x})\ev(d_0,\vec{x}) \neq 0$ \qquad $(\bullet)$ \\
& \qquad $\land\,(k)\,\neg\,\Pro_\T
                  (k,\code{(\vec{x})\,\mr{ev}(d_0,\vec{x}) \neq 0})$
}
{ $\T \derives\ \Pro_\T
     (\num\,\ulj,\code{(\vec{x})\,\mr{ev}(d_0,\vec{x}) \neq 0})$ \\
& $\land\,(k)\,\neg\,\Pro_\T
                (k,\code{(\vec{x})\,D(\vec{x}) \neq_Z 0})$
}
a contradiction: 
appropriate $\ulj$ is available from $(\bullet)$ via  
derivation-to-$\Proof$-internalisation (\emph{g\"odelisation}).

\bigskip
$[\,$For the time being we consider $\T$ as frame, not (yet) 
$\T^+=\T+\Con_\T.\,]$

\bigskip
For $\T = \bfQ$ quantified, with (countable) \emph{axiom of choice} 
$\ACC,$ in particular $\bfQ = \PA+\ACC$ Peano Arithmetic with choice,
we define the \emph{undecided part} of $\mr{DIO}$ as
\begin{eqnarray*}
\Psi &=& \Psi^{\bfQ} \\    
&=& \set{d \in \mr{DIO}:\forall\,k\ \mr{ev}(d,\ct_*\,k) \neq 0 \\
&& \,\land\,\forall\,k
  \ \neg\,\Pro_{\bfQ}(k,\code{(\vec{x})\,\mr{ev}(d,\vec{x}) \neq 0)})} \\
&\subset& \mr{DIO} = \Zlist \subset \N. 
\end{eqnarray*}
With this definition we get 
\begin{align*}
\bfQ \derives\
& \Psi \neq \emptyset
    \implies \choice_{\Psi}: \one \to \Psi \subset \N\ \emph{total} \\
& (\text{choice available by $\ACC:$ non-empty sets have 
        \emph{defined points}}) \\
& \implies \mu\set{d:t(d)\ \text{non-terminating}}:\one \to \Psi\ \emph{total.}  
\end{align*}
This means: the {assumption} of (formal) \emph{existence} of a 
$d \in \mr{DIO}$ for which decision race $t: \mr{DIO} \parto \N$ does \emph{not} terminate, leads to a (\emph{defined}) point 
  $$d_0:\one \to \mr{DIO}$$ 
for which $t$ derivably does not terminate.

But this is \textbf{excluded} by pointwise non-derivability above 
of non-termination, within frame $\bfQ$ assumed consistent. 

So we have shown
\begin{align*}
& \bfQ,\PA+\ACC \derives\ \Psi = \emptyset,\ \text{i.\,e.} \\ 
& \bfQ \derives\ (\forall d\in \mr{DIO})[\exists k\,\ev(d,\mr{ct}_*k)=0 \\
& \qquad\qquad
  \lor \exists k\,\Pro_{\DIO}(k,\code{(\vec{x})\,\ev(d,\vec{x}}) \neq 0)],
\end{align*}
whence

\smallskip
\textbf{Termination Theorem:} $\bfQ,\ZFC,\PA+\ACC$ derive race $t$ to 
terminate on all diophantine codes $d,$ on all $d \in \mr{DIO} = \Zlist.$

\section{Correct termination of decision $\nabla$}

\textbf{In particular} ($\bfQ^+ = \bfQ+\ACC$ stronger than $\bfQ$): 
\begin{align*}
& \bfQ^+\ \textbf{derives} \\
& \quad 
    \text{overall termination of $\mu$-recursive} \\ 
& \quad
    \text{termination race}\ 
                         t=t^{\bfQ}(d):\mr{DIO} \to \N: \\
& \bfQ^+\ \derives\
    [\,(\forall\,d \in \mr{DIO})\ t(d) \in \N
                   \ \text{{\emph{defined}}}\,]
\end{align*}

\textbf{Hence,} by Decision Correctness within $\bfQ^+:$

\smallskip
$\bfQ^+$ {$\mathbf{derives}$} \\ 
{overall \emph{correct} termination} of $\mu$-recursive \emph{decision}\\ 
$\nabla: \mr{DIO} \to \two,$ \textbf{main result} here:
\begin{align*}
& \nabla(d) \\
& = \begin{cases}
      0\ \myif\ \mr{ev}(d,t(d)) = 0 \\
      \qquad 
        \ [ \ \implies d\ \emph{has}
                         \ \text{a zero}\ \vec{z} \in \Z^*\ ] \\
      1\ \myif\ \Pro_{\DIO} 
              (t,\code{(\forall\,\vec{x})
                      \ \mr{ev}(d,\vec{x}) \neq 0}) \\
      \qquad 
        \ [ \ \implies d\ \text{has}\ \emph{no}\ \text{zero}\ ] 
    \end{cases}: \mr{DIO} \to \two.
\end{align*}

\section{Comparison with Matiyasevich's\\negative result}

\emph{Main result} above says in terms of the theory $\mathbf{TM}$ 
of TURING machines, by the established part of CHURCH's thesis:

\smallskip
\emph{For concrete diophantine polynomials} $D = D(\vec{x}): \Z^m \to \Z:$ 

For quantified arithmetical choice theories $\bfQ+\ACC$ like 
$\ZFC$ and already $\PA+\ACC,$ 

$\bfQ^+ = \bfQ+\Con_{\bfQ}$ {$\mathbf{derives}$}:

\smallskip
\emph{TURING machine ${{\mrTM}}_{\nabla_{\bfQ}}$ 
corresponding---CHURCH---to 
totally defined $\mu$-recursive decision map 
  $$\nabla_{\bfQ}: \mr{DIO} \to \set{0,1},$$ 
\emph{when written} \emph{coefficient list} $\ccode{D}$ 
of a diophantine polynomial $D$ on its (initial) TAPE,  
eventually \textbf{reaches} \emph{HALT state,}
leaves result $0$ (as its \emph{final} TAPE) 
$\mathbf{iff}$ $D$ \emph{has} a zero\\
$\vec{z}:$ $D(\vec{z}) = 0,$  
and \textbf{result} $1$ $\mathbf{iff}$ $D$ is overall \emph{non-null}: \\ 
$(\forall\,\vec{x} \in \Z^*)\,[\,D(\vec{x}) \neq 0\,].$}

\medskip
This contradicts \textbf{Matiyasevich's} THEOREM \emph{unsolving} 
Hilbert's 10th problem, within theory $\bfQ^+$ which strengthens
his framework of Peano Arithmetic $\PA+\ACC$ with countable axiom
of choice. Whence

\subsubsection*{Conclusion:}
\begin{itemize}
\item
$\ZFC^+ = \ZFC+\ConZFC$ is contradictory, so

\item
$\ZFC\derives\neg\,\ConZFC:$ $\ZFC$ \emph{is internally inconsistent,}

\item
same for theory $\PA+\ACC:$ 
 
\emph{Peano-Arithmetic with axiom of countable choice is internally
inconsistent}

\item
\textbf{Question:} is already Peano Arithmetic $\PA$ by itself 
internally inconsistent? It would be if axiom $\ACC$ of countable 
choice were derivable within $\PA$ or independent from $\PA,$ 
as is axiom of choice $\AC$ from \textbf{set} theory. This would
mean that formal existential quantification is incompatible with
free-variables Primitive Recursive Arithmetic $\PR.$
\end{itemize}

\subsubsection*{Discussion}
\begin{itemize}
\item 
After his talk at Humboldt University Berlin, I have mailed to
Matiyasevich the question, if his \emph{unsolving} of Hilbert's
10th problem is really constructive: it depends heavily on
formal existential quantification. No reply: may be he considers this
question when present paper will be brought to his attention.
 
\item
I have submitted the 200? version of present work, claiming 
self-inconsistency $\PA\derives \neg\,\Con_{\PA},$ to the 
\emph{Journal of Symbolic Logic.} The (anonymous) referee: 

\emph{... this is certainly false. ...} Robert 'Rob' \NAME{Goldblatt} ed.:
\emph{under these circumstances etc.}
\end{itemize}

\textbf{What is such editorial policy good for?}

\section{Hilbert 10 constructively}

In this section we show that the \emph{local} version $\nabla[D]: 1 \parto 2$ 
of the $\mu$-recursive \emph{decision algorithm} 
$\nabla = \nabla_{\DIO}(d): \mathit{DIO} \parto \two$ 
\emph{irrefutably} \emph{decides} \emph{each (single)} diophantine 
equation---\emph{correctly}---when placed in p.\,r.
\emph{non-infinite-descent theory} $\piR=\PR+(\pi)$ of
op.\,cit.\  in the References. 


This will give a positive solution to Hilbert's 10th problem
in that constructive framework, at least when
stated in its original form quoted in first section above.



\medskip
Formally, this \textbf{problem} allows for solution by a separate decision algorithm  (``process'') for each diophantine polynomial.
By \emph{localisation} at a given polynomial, we extract such a 
decision-\emph{family} from the forgoing sections, and formalise it 
within $\piR.$

We index that family (externally) by the 
\emph{diophantine constants} $\delta:\one \to \mr{DIO} \subset \N,$ 
among which the diophantine polynomials 
$$D = D(\vec{x}) = D(x_1, \ldots, x_{\bs{m}}): \Z^{\bs{m}} \to \Z$$ 
are represented by their coefficient list codes 
$\ccode{D}: \one \to \mr{DIO}.$ 

\medskip
\textbf{Definition:} For PR predicates $\ph_0, \ph_1: A \times \N \to \two$
we define the \emph{race winner predicate} 
 $$\mu_{\lor}[\ph_0,\ph_1]: A \to \two$$ 
between $\ph_0$ and $\ph_1$ slightly assymmetrically by
\begin{eqnarray*}
&& \mu_{\lor} [\ph_0,\ph_1] = \mu_{\lor} [\ph_0,\ph_1](a)\\
&& \defeq
(\mathit{dc} \circ (\ph_0,\ph_1)) 
  \parcirc (A \times \mu[\ph_0 \,\lor\, \ph_1]) \parcirc \Delta_A: \\
&& A \to A \times A \parto A \times \N \to \two \times \two 
       \overset{\mathit{dc}} {\longrightarrow} \two, \ \text{with} \\
&& \mathit{dc} = \mathit{dc}(u,v): \two \times \two \parto \two 
   \ \text{defined by} \\
&& \mathit{dc}(u,v) \defeq
\begin{cases}
0 \ \text{if}\ u = 1, \\
1 \ \text{if}\ u = 0\ \land\ v = 1, \\
\text{\emph{definably undefined} if}\ u = v = 0.
\end{cases}
\end{eqnarray*}
This (partial) race winner predicate 
$\mu_{\lor}[\ph_0,\ph_1](a): A \parto \two$ 
is characterised---within $\bs{S}=\PR$ as well as in $\bs{S}=\piR$---by
\begin{align*}
  \bs{S} \derives\, 
           & [\,\ph_0 (a,n) 
                \,\land\, \underset{i<n} {\land}\,\neg\,\ph_1(a,n)
                 \implies \mu_{\lor}[\,\ph_0,\ph_1\,](a)=0\,] \\
           & \land\, [\,\ph_1(a,n)
             \,\land\,\underset{i \leq n}{\land}\,\neg\,\ph_0(a,n) 
                \implies \mu_{\lor}[\,\ph_0,\ph_1\,](a)=1\,]. 
\end{align*}
We allow us to write for this intuitively---in classical terms of a (partial)
case-distinction: 
$$
  \mu_{\lor}[\,\ph_0,\ph_1\,](a) = 
  \begin{cases}
    0 \ \text{if}\  \mu\ph_0(a) < \infty
                       \,\land\,\mu\ph_0(a) \leq \mu\ph_1(a), \\ 
    1 \ \text{if}\  \mu\ph_1(a) < \infty 
                             \,\land\,\mu\ph_1(a) < \mu\ph_0(a).
  \end{cases}
$$
Our decision family 
  $$\nabla[\delta]: 1 \parto \two,\ \delta: \one \to \mr{DIO} \subset \N$$ 
now is defined in the present $\mu$-recursive
frame as this type of race winning, of PR search for a zero 
(in the evaluation) of $\delta$ against PR search for a (first) 
internal non-nullity \emph{proof} for (the evaluation) of $\delta,$ 
namely by
\begin{align*}
\nabla[\delta] 
& \defeq \mu_{\lor} [\ph_0[\delta],\ph_1[\delta]]: 1 \parto \two, 
  \ \text{with} \\
\ph_0[\delta](k) 
& \defeq [\,\ev(\delta,\ct_*(k)) = 0\,]: \N \to \two, \\ 
\ph_1[\delta](k) 
& \defeq \Pro_{\bs{S}}(k,\code{(\vec{x})\ev(\delta,\vec{x}) \neq 0}.
\end{align*}
Here
  $$\ev = \ev(d,x): \N \times \N \supset \mr{DIO} \times \Z^* \to \Z$$
is evaluation with the characteristic \textbf{evaluation property}
\begin{align*}
& \ev(\ccode{D},(x_1,\ldots,x_{\bs{m}}))
  = D(x_1, \ldots, x_{\bs{m}}): Z^{\bs{m}} \to \Z,
\end{align*} 
realised by (iterated) \NAME{Horner}'s schema (each application reduces 
the number of remaining variables by 1), or by ``brute force'' evaluation of monomials.

\subsection{Decision Correctness}

\textbf{Soundness Recall:} Main result of op.\,cit.\  
in the References is (logical) \emph{soundness} of theory $\piR:$
\begin{itemize}
\item
For a (p.\,r.\,) predicate  $\chi = \chi(a): A \to \two$ we have
$$\piR\,\derives\ \Pro_{\piR}(k,\code{\chi})  
  \implies \chi(a): \N \times A \to \two,$$ 
$a \in A$ free, meaning here \emph{for all} $a\in A,$ and $k \in \N$ free, 
meaning here \emph{exists} $k\in\N.$ This entails 

\item
\emph{$\PR$ soundness} of $\piR:$ For a p.\,r.\ predicate  
$\chi = \chi(a): A \to \two,$
$$\piR\,\derives\ \Pro_{\PR}(k,\code{\chi})  
  \implies \chi(a): \N \times A \to \two,$$ 
as well as in particular

\item
\emph{Diophantine soundness} of $\piR:$
for a diophantine polynomial $D=D(\vec{x}):\Z^*\to \two$
$$\piR\,\derives\ \Pro_{\piR}(k,\code{(\vec{x})D(\vec{x})\neq 0})  
  \implies D(\vec{x})\neq 0,$$ 
$k\in \N,\ \vec{x}\in \Z^*$ free.

\item
Already $\PR^+ = \PR+\Con_{\PR}$ is diophantine sound. 
This needs an extra Proof. 

We consider here frame $\bs{S}=\piR,$
$$\piR^+ = \piR+\Con_{\piR} = \piR,$$ 
the latter by op.\,cit.\ equivalent to soundness of theory $\piR.$                                  
\end{itemize}

Namely from PR Soundness we get the

\bigskip
\textbf{Local Correctness-Lemma} for $\nabla[\delta]$ in $\piR:$ 
The partial $\PR$-map $\nabla[\delta]:\one \parto \two$
has the following correctness properties:

\smallskip
$\piR \derives\ :$
\begin{itemize}
\item $\delta$ does not fall in 
\emph{both} of the two defined-cases stated for $\nabla[\delta],$

\item 
$\nabla[\delta] = 0 
 \implies \ev(\delta,\ct_{*}\circ \mu\ph_0 [\delta])=0:$  
$\delta$ is implied to have available a zero in its \emph{evaluation,}
  
\item 
$\nabla[\delta] = 1  \implies \ev(\delta,\vec{x}) \neq_Z 0,$ 
$\vec{x}$ free in $\Z^*$: 
$\delta$ is implied to be evaluated globally non-null, in particular:

\item By diophantine evaluation for 
$D = D(x_1, \dots x_{\bs{m}}): \Z^* \to \Z$ diophantine: 
\begin{itemize}
\item 
$\nabla[D] := \nabla[\ccode{D}] = 0 
 \implies D(\ct_*(\mu\ph_0 [\ccode{D}]))=0:$ 
 
$D$ is implied to have a zero, as well as 
  
\item $\nabla[D]=1 \implies [\,D(\vec{x}) \neq 0\,],$ 
here again $\vec{x}$ free over $\Z^*:$

$D$ is implied to be globally non-null\ \,\textbf{\qed}
\end{itemize}   
\end{itemize}



\subsection{Decision Termination}

The final question to treat for this---canonical---family
$$\nabla = \nabla_{\mbf{DIO}}[\delta]:\one \parto \two, 
  \ \delta: \one \to \mr{DIO} \subset \N$$
of \emph{local}---$\mu$-recursive---decision algorithms, is
\emph{termination,} for each $\delta,$ in particular for
$\delta = \ccode{D},$ $D = D(\vec{x})$ diophantine. 

\smallskip
\textbf{Assume} $\nabla[d_0]$ \emph{not} to terminate for a particular 
\emph{constant} $d_0: \one \to \mr{DIO},$ in particular $d_0$ of
form $D_0 = D_0(\vec{x}).$  

Since we argue here purely \emph{syntactically}---within the \emph{theory} 
$\widehat{\bs{S}} \bs{\supset} \bs{S} \bs{=} \PR+(\mr{abstr})$ of \emph{partial} 
p.\,r.\ maps---no modelling in mind except some primitive recursive 
\emph{Meta}mathematics (these in turn g\"odelised within $\bs{S}$)---we 
discuss the stronger assumption

$\nabla[d_0]$ $\T$-\emph{derivably} 
does \emph{not} terminate for a given diophantine constant 
$d_0: \one \to \mathit{DIO},$ $\T$ an extension of $\bs{S}.$

\medskip
This \textbf{assumption} reads:
$$\T \derives\,(k)\psi[d_0](k):$$
here $k$ is free over $\N,$ and the PR predicate
$\psi[d_0](k): \N \to \two$ is defined by
\begin{align*}
& \psi[d_0](k) = \psi_0[d_0](k) \land \psi_1[d_0](k)\ \text{with} \\
& \psi_0[d_0](k) = [\,\ev(d_0,\ct_*(k)) \neq 0 \,], \text{and} \\
& \psi_1[d_0](k) = \neg\,\Pro_{\T}(k,\code{\ev(d_0,\vec{x})\neq 0}).
\end{align*}
So the assumption (``of the contrary'') reads:
\begin{align*}
\T \derives\, 
& [\,\ev(d_0,\ct_{*}(k) ) \neq 0 \,] \\
& \land \neg\,\Pro_{\T}(k,\code{(\vec{x})\ev(d_0,\vec{x}) \neq 0}).
\end{align*}
Here $k\in\N$ is the only free variable in the \emph{accessible} level, 
$\vec{x}$ is free over $\Z^*,$ but \emph{encapsulated} within g\"odelisation,
\emph{not visible} on the object language level.

\smallskip
The derivably-non-termination assumption 
$$\T \derives\, \psi[d_0](k),\ k\ \free,$$
would entail in particular (first conjunct $\psi_0[d_0]$):
$$\T \derives\, \ev(d_0,\ct_*(k)) \neq 0:\N\to\two.$$

\emph{Internalising} (\emph{formalising}) this metamathematical statement, 
we (would) get by Proof-Internalisation---\cf \NAME{Smorynski} 1977---a 
\emph{constant} $p_0: \one \to \mathit{Proof}_\T \subset \N$ 
\emph{guilty} for this last statement:
$$\T \derives\ \Pro_\T(p_0,\code{\ev(d_0,\vec{x}) \neq 0});$$
this would give, by definition of $\nabla[d_0]:$ 
$$\T \derives \nabla[d_0] = 1,$$
a contradiction to our assumption that $d_0$ be derivably 
\emph{not decided} by $\nabla_{\DIO},$ \ie to $\T \derives \psi[d_0].$

\bigskip
\textbf{Conclusion:}

\begin{itemize}
\item $\piR=\piR+\Con_{\piR}$ derives the 
alleged decision algorithm (family)  
$\nabla = \nabla_{\DIO}[D]:\one \parto \two$ to be \emph{correct}
for each diophantine polynomial (if defined).

\item no diophantine polynomial $D=D(\vec{x})$ 
can come with a $\T$-proof (i.\,p.\ a $\piR$-proof) showing $\nabla[D]$ 
to be \emph{undefined,} \emph{not} to terminate, 
in other words: 

\item \emph{correct termination} of the 
$\mu$-recursive \emph{decision family} $\nabla = \nabla_{\DIO}[D]$ 
at each diophantine polynomial is $\piR$-\emph{irrefutable,} 
in the sense that \textbf{otherwise}---refutation---
\begin{align*}
& \piR \derives \Pro_{\piR}(q,\code{\mr{false}}),  
  \ q: \one \to \N\ \text{a suitable PR point,}                                                                                    
\end{align*}
inconsistency of (self-consistent) theory $\piR$ would be the consequence. 

\end{itemize}





\subsection*{Outlook}

Irrefutable correct termination of \emph{uniform} decision algorithm
$$\nabla_{\DIO}=\nabla_{\DIO}(d):\mr{DIO}\parto\two,\ d\in\mr{DIO}\ \free$$
is treated within the general framework of 

\emph{Arithmetical Decision} 
to come.

\end{document}